\documentclass[final]{svjour3PP}

\usepackage{amsfonts,amsmath,amssymb}

\usepackage{mathptmx}
\usepackage{cite}

\AtBeginDocument{%
  \setlength{\oddsidemargin}{\dimexpr(\paperwidth-\textwidth)/2-1in}%
  \setlength{\evensidemargin}{\oddsidemargin}%
  \setlength{\topmargin}{%
    \dimexpr(\paperheight-\textheight)/2-\headheight-\headsep-1in}%
}

\def\0{\ensuremath{\mathbf{0}}}

\usepackage{graphicx}
\usepackage{a4,amssymb,amsmath}

\def\be{\begin{equation}}
\def\bea{\begin{eqnarray*}}
\def\ee{\end{equation}}
\def\eea{\end{eqnarray*}}
\def\ba{\begin{array}}
\def\ea{\end{array}}
\def\bi{\begin{itemize}}
\def\ei{\end{itemize}}
\newtheorem{theo}{Theorem}
\newtheorem{lem}{Lemma}

\title{Stochastic Subspace Correction in Hilbert Space}
\author{Michael Griebel \and Peter Oswald}
\institute{M. Griebel\at Institute for Numerical Simulation, Universit\"at Bonn, Wegelerstr. 6, 53115 Bonn, and 
Fraunhofer Institute for Algorithms and Scientific Computing (SCAI), Schloss Birlinghoven, 53754 Sankt Augustin\\
Corresponding author, tel.: +49-228-733437, fax: +49-228-737527,\\ 
\email{griebel@ins.uni-bonn.de}
\and P. Oswald \at Institute for Numerical Simulation, Universit\"at Bonn, Wegelerstr. 6, 53115 Bonn,\\
\email{agp.oswald@gmail.com}}

\titlerunning{Stochastic Subspace Correction in Hilbert space}
\authorrunning{M. Griebel \and P. Oswald}

\date{}
\begin{document}
\maketitle
\begin{abstract}
We consider an incremental approximation method for solving variational problems in infinite-dimensional separable Hilbert spaces, where in each step a randomly and independently selected subproblem from an infinite collection of subproblems is solved. We show that convergence rates for the expectation of the squared error can be guaranteed under weaker conditions than previously established in \cite{GrOs2016}.

\keywords{infinite space splitting \and subspace correction \and multiplicative Schwarz \and block coordinate descent \and    greedy \and randomized \and convergence rates \and online learning}
\subclass{65F10 \and 65N22\and 49M27}
\end{abstract}

\section{Introduction}\label{sec1}
The fast solution of quadratic minimization problems or, correspondingly, of large linear systems of equations is an important topic in many application areas of numerical simulation. To this end, iterative algorithms play a major role. They can be formalized by means of subspace correction methods, either in the so-called additive or the multiplicative variant, see \cite{Xu,GrOs1995}. There, a set of variational problems in appropriate subspaces is chosen and the current approximation is iteratively (either collectively or successively) improved by means of the respective solutions of these subproblems. The subspaces can be one-dimensional or, in a block type fashion, they can have arbitrary (finite) dimension as well. Examples are the well-known Jacobi and Gauss-Seidel algoritms from linear algebra and their block-wise variants, but also domain decomposition methods or multigrid and multilevel techniques from scientific computing.

For a (finite or infinite) set of subproblems at hand, the question arises in which order the incremental updates should be made and what the convergence behavior of the associated subspace correction algorithm will be. Most conventionally, the order is a priorily fixed in a deterministic fashion. This is the case for basically all the classical methods, like Gauss-Seidel, domain decomposition or multigrid methods. The order of traversal through the subproblems is prescribed by the method itself. Examples are lexicographical or so-called red-black orderings  for systems stemming from finite element or finite difference discreitizations and, additionally, level wise traversal orderings in multigrid algorithms. Besides, for the multiplicative subspace correction approach, there are greedy methods where the next subspace is identified according to an optimization criterion such that the actual error is reduced by the following incremental update as much as possible. This may substantially improve the convergence of the overall algorithm. A detailed analysis of various greedy approximation methods is given in the seminal book \cite{Tem2012}. A simple example from linear algebra is the so-called Gauss-Southwell approach, where the next update variable is that with the largest residuum. Usually, in the case of finitely many subspaces, the determination of the optimal next subspace can be done exactly, but it involves additional costs. In the case of infinitely many subspaces, this is not possible any more, and heuristic choices are employed there in practical methods.

Besides a deterministic or greedy pick, we may also choose the next subproblem in a random fashion according to a 
probability distribution $\rho$ on the set of subspaces, see \cite{GrOs2012} and the references cited therein. 
The analysis of such stochastic iterations has been a very active research topic in large-scale convex optimization, see \cite{FerRic2016} for a recent survey, but also in the area of machine learning and compressed sensing.
Compared to the greedy approach, the cost for determining the next subspace is dramatically reduced to the cost of sampling the underlying probability distribution $\rho$. Moreover, the random pick is feasible also for infinite sets of subspaces.  But the question is now what the associated convergence rate (in expectation) will be. For finitely many subspaces, the answer is very encouraging \cite{GrOs2012}: Both greedy and stochastic iterations yield the same exponential rates of convergence, although with different constants, and in the latter case almost surely and in expectation only.

In this article, we deal with the case of an infinite number of subspaces for which a first comparison of greedy and stochastic subspace correction methods was carried out in \cite{GrOs2016} for countable sets of subspaces. It was shown that the (much more involved and costly) greedy method converges at an algebraic rate for solutions from a certain class ${\mathcal{A}}_1$ while  basically the same convergence rate can be achieved in expectation by a stochastic subspace correction method on a class ${\mathcal{A}}_\infty^\rho\subset {\mathcal{A}}_1$ depending on $\rho$. Details will be given in the next sections.

The aim of this paper is to show that convergence rates for the expectation of the squared error can be guaranteed under weaker conditions than previously established in \cite{GrOs2016}, namely for solutions from a class $A_2$ still depending on $\rho$, where  ${\mathcal{A}}_\infty^\rho \subset A_2 \subset {\mathcal{A}}_1$. This result reveals some connection to the theory of approximation algorithms in reproducing kernel Hilbert spaces, and may also allow for a wider range of applications of incremental, multiplicative subspace correction methods with randomly picked orderings which may have interesting applications in numerical linear algebra, scientific computing, quadratic optimization, machine learning and compressed sensing.

The remainder of this paper is organized as follows:
In Section \ref{sec2} we give basic notation and introduce our multiplicative subspace correction/approximation algorithm with random picking in the case of a family of subspaces
$\{V_\omega:\,\omega\in \Omega\}$ with an infinite (possibly uncountable) index set $\Omega$. Moreover, we give in Theorem \ref{theo1} and Theorem \ref{theo2} sharp bounds of its error and thus of its convergence rate in expectation for the class $A_2$.
In Section \ref{sec3} we discuss various examples of our abstract theory. First, we consider the case of a countable index set $\Omega$ and discrete probability measures $\rho$ on it. Moreover, in Lemma \ref{lem1} we also relate our new function class $A_2$ to the classes ${\mathcal A}_\infty^\rho$ and ${\mathcal{A}}_1$, previously used in \cite{GrOs2016}.  Then, we consider the case of stochastic approximation in reproducing kernel Hilbert spaces and show that our theory can be applied there as well. Next, we study the case of general unit norm dictionaries and approximation with these, and provide in Theorem \ref{theo3} a version of our main results from Section \ref{sec2} with simplified proof.
Finally, we deal with a collective approximation problem from \cite{AJOP2017} and show how our theory applies. 
We conclude in Section \ref{sec4} with some further remarks on  our convergence results.

\section{Details and Proofs}\label{sec2}

Throughout this paper, let $V$ be a separable real Hilbert space.  For a given continuous and coercive Hermitian form $a(\cdot,\cdot)$ on $V$ and a bounded linear functional $F$ on $V$, we consider the variational problem of finding the unique element $u\in V$ such that
\be\label{VP}
a(u,v)=F(v)\qquad\forall v\in V.
\ee
Equivalently, (\ref{VP}) can be formulated as quadratic minimization problem in $V$ or as linear operator equation in the dual space of $V$. In the following, we use the fact that 
$a(\cdot,\cdot)$ defines a spectrally equivalent scalar product on $V$, equip $V$ with it, and write $\|v\|={a(v,v)}^{1/2}$. 

Our aim is to study a particular instance of an incremental subspace correction (or Schwarz iterative) method for solving (\ref{VP}). Let $\Omega$ be a fixed index set equipped with a probability measure $\rho$ 
(compared to \cite{GrOs2016}, we also allow for uncountable $\Omega$, see below for an example). Consider 
a family $\{V_\omega\}_{\omega\in\Omega}$ of separable real Hilbert spaces, each equipped with a spectrally equivalent scalar product 
$a_\omega(\cdot,\cdot)$ and norm $\|v_\omega\|_\omega:=a_\omega(v_\omega,v_\omega)^{1/2}$, and linear operators $R_\omega:\,V_\omega\to V$ such that
\be\label{La}
\|R_\omega\|_{V_\omega\to V} = \sup_{\|v_\omega\|_\omega=1} \|R_\omega v_\omega\| \le \Lambda < \infty,\qquad \omega\in \Omega.
\ee
Finally, we introduce another family 
of linear operators 
$T_\omega:\,V\to V_\omega$ by the solution of auxiliary variational problems in $V_\omega$:
\be\label{VPo}
a_\omega(T_\omega v,v_\omega)=a(v,R_\omega v_\omega)\qquad \forall\;v_\omega\in V_\omega,\qquad \omega\in \Omega.
\ee
It is easy to see that $\|T_\omega\|_{V\to V_\omega} \le \Lambda$ as well. Without loss of generality, we can assume
that $\mathrm{Ker}(R_\omega)=\{0\}$ for all $\omega$ (otherwise replace $V_\omega$ by $V_\omega \ominus_\omega
\mathrm{Ker}(R_\omega)$).

With these preparations at hand, the algorithm under consideration has the form
\be \label{Rec}
u^{(m+1)} = \alpha_m u^{(m)} + \xi_m R_{\omega_m}r^{(m)}_{\omega_m},\qquad r^{(m)}_{\omega_m}=T_{\omega_m}(u-u^{(m)}),\quad m=0,1,\ldots,\quad u^{(0)}=0,
\ee
where $\{\omega_m\}$ is a sequence of independent samples from $\Omega$ which are identically distributed according to $\rho$. Furthermore $\alpha_m=1-(m+2)^{-1}$, and $\xi_m$ is such that the error 
$$
\delta_{m+1}^2:=\|u-u^{(m+1)}\|^2
$$ 
is minimized. This gives the explicit formula
\be\label{Xi}
\xi_m=\mathrm{argmin}_\xi \|u-\alpha_m u^{(m)}-\xi R_{\omega_m}r^{(m)}_{\omega_m}\|^2 =\frac{F(R_{\omega_m}r^{(m)}_{\omega_m})-\alpha_m a(u^{(m)},R_{\omega_m}r^{(m)}_{\omega_m})}{a(R_{\omega_m}r^{(m)}_{\omega_m},R_{\omega_m}r^{(m)}_{\omega_m})}.
\ee
Since $r^{(m)}_{\omega_m}$ is defined via (\ref{VP}) and (\ref{VPo}) by the variational problem
\be\label{VPoe}
a_\omega(r^{(m)}_{\omega_m} ,v_{\omega_m})=a(u-u^{(m)},R_{\omega_m} v_{\omega_m})=F(R_{\omega_m} v_{\omega_m})-a(u^{(m)},R_{\omega_m} v_{\omega_m})\qquad \forall\;v_\omega\in V_\omega,
\ee
we see that (\ref{Rec}) can be executed once $u^{(m)}$ and $\omega_m$ are available.

We note that this way $u^{(m)}$ and thus $\delta_m^2$ become random variables on $\Omega^m$ equipped with the product measure $\rho^m$. To provide estimates for the expected squared error $\mathbb{E}(\delta_m^2)$, we need the notion of Bochner integrals \cite{Bo}. Given any Bochner-measurable $V$-valued
function $\phi:\, \omega\in \Omega \to \phi_\omega \in V$, its Bochner integral
\be\label{E}
\mathbb{E}_\rho (\phi):=\int_{\Omega} \phi_\omega\,d\rho_\omega
\ee
is well-defined with value in $V$ iff the scalar integral
\be\label{En}
\mathbb{E}_\rho (\|\phi\|):=\int_{\Omega} \|\phi_\omega\| \,d\rho_\omega < \infty
\ee
exists. The Bochner integral is similarly well-defined if $V$ is replaced by a separable Banach space. In the case of a discrete probability measure on a countable index set $\Omega$, measurability of $\phi$ is not an issue, in other situations, it needs to be checked.
For the following, we assume that for any fixed $e\in V$ the function 
\be\label{A1}
\tilde{\psi}:\, \omega\in \Omega \to \tilde{\psi}_\omega \in 
R_{\omega}(V_\omega)\subset V,\qquad \tilde{\psi}_\omega := \left\{
\ba{ll}\frac{R_{\omega}T_\omega e}{\|R_{\omega}T_\omega e\|},& R_{\omega}T_\omega e\neq 0,\\
0,& R_{\omega}T_\omega e = 0,\ea\right.
\ee
is Bochner-measurable. 

Next, we introduce the class $A_2\equiv A_{2,\rho}\subset V$ which will play a central role in the convergence theory for (\ref{Rec}). We say that $u\in V$ belongs to $A_2$ if there exists 
a Bochner-measurable function $\phi: \,\omega\to R_\omega v_\omega$ with $v_\omega\in V_\omega$ for all $\omega\in\Omega$ such that the scalar-valued function $\omega\to \|v_\omega\|_\omega$ is also measurable, and
\be\label{A2}
u=\mathbb{E}_\rho(\phi)=\int_{\Omega} R_\omega v_\omega\,d\rho_\omega,\qquad \mathbb{E}_\rho(\|v_\omega\|_\omega^2)=\int_{\Omega} \|v_\omega\|_\omega^2\,d\rho_\omega<\infty,
\ee
Define a norm  on $A_2$ by
\be\label{NA2}
\|u\|_{{A}_2}:= \inf\, \mathbb{E}_\rho(\|v_\omega\|_\omega^2)^{1/2},
\ee
where the infimum is taken with respect to all admissible representations of $u$ in (\ref{A2}).
How this class is related to the classes $\mathcal{A}_p^\gamma$ introduced in \cite{GrOs2016} for discrete measures
$\rho$ on countable index sets $\Omega$ and other classes
used in similar context in the literature will be elaborated on in Section \ref{sec3}.

The central result of this note is the following:
\begin{theo}\label{theo1}
If (\ref{A1}) holds and if $u$ belongs to the linear space ${A}_2$ induced by the condition (\ref{A2}) 
then, for the incremental approximation algorithm (\ref{Rec}), we have 
\be\label{EC2}
\mathbb{E}(\delta_{m}^2) \le \frac{(\Lambda\|u\|_{{A}_2}+\|u\|)^2}{m+1},\qquad m=0,1,\ldots.
\ee
\end{theo}

{\bf Proof }. We start with an analysis of the error reduction in one recursion step, i.e., with an estimate
of $\mathbb{E}_\rho (\delta_{m+1}^2|u^{(m)})$. By (\ref{Xi}) and with the notation 
$$
e^{(m)}:=u-u^{(m)}, \qquad \bar{\alpha}_m:=1-\alpha_m=(m+2)^{-1}, \qquad 
w:=\alpha_m e^{(m)}+\bar{\alpha}_m u,
$$ 
we have
\bea
\delta_{m+1}^2&=&\min_\xi\, \|\alpha_m(u-u^{(m)})+\bar{\alpha}_m u - \xi R_{\omega_m}r^{(m)}_{\omega_m}\|^2\\
&=& \|w\|^2 - a(w,\tilde{\psi}_{\omega_m})^2 = \alpha_m^2 (\delta_m^2-a(e^{(m)},\tilde{\psi}_{\omega_m})^2)\\
&&\quad +2\alpha_m\bar{\alpha}_m(a(e^{(m)},u)-a(e^{(m)},\tilde{\psi}_{\omega_m})a(u,\tilde{\psi}_{\omega_m}))
+\bar{\alpha}_m^2(\|u\|^2-a(u,\tilde{\psi}_{\omega_m})^2).
\eea
Here and throughout the proof $\tilde{\psi}$ stands for the function defined in
(\ref{A1}) for $e=e^{(m)}$. The measurability assumption for this $\tilde{\psi}$ allows us to take expectations with respect to the choice of $\omega_m$ in the above error representation:
\be\label{I1}
\ba{l}
\mathbb{E}_\rho(\delta_{m+1}^2|u^{(m)})= \alpha_m^2 (\delta_m^2-\mathbb{E}_\rho(a(e^{(m)},\tilde{\psi}_{\omega})^2))\\
\\
\qquad +2\alpha_m\bar{\alpha}_m(a(e^{(m)},u)-\mathbb{E}_\rho(a(e^{(m)},\tilde{\psi}_{\omega})a(u,\tilde{\psi}_{\omega})))
+\bar{\alpha}_m^2(\|u\|^2-\mathbb{E}_\rho(a(u,\tilde{\psi}_{\omega})^2)).
\ea
\ee
By the definition of $r_\omega^{(m)}=T_\omega  e^{(m)}$ via (\ref{VPo}) we have
\be\label{I0}
\frac{\|R_\omega r_\omega^{(m)}\|}{\|r_\omega^{(m)}\|_\omega} a(e^{(m)},\tilde{\psi}_\omega) = a_\omega(r_\omega^{(m)},\frac{r_\omega^{(m)}}{\|r_\omega^{(m)}\|_\omega})\ge a_\omega(r_\omega^{(m)},\frac{v_\omega}{\|v_\omega\|_\omega}) = \frac{a(e^{(m)},R_\omega v_\omega)}{\|v_\omega\|_\omega}
\ee
for any $v_\omega\in V_\omega$ and $\omega\in \Omega$. Together with (\ref{La}) and (\ref{A2}), this implies
$$
a(e^{(m)},u) = \mathbb{E}_\rho(a(e^{(m)},R_\omega v_\omega)) \le \Lambda \mathbb{E}_\rho(\|v_\omega\|_\omega a(e^{(m)},\tilde{\psi}_{\omega})).
$$

Thus, we can apply the Cauchy-Schwarz inequality to the second term in the right-hand side of (\ref{I1}):
\bea
&& 2\alpha_m\bar{\alpha}_m( a(e^{(m)},u)-\mathbb{E}_\rho(a(e^{(m)},\tilde{\psi}_{\omega})\,a(u,\tilde{\psi}_{\omega}))\\
&&\qquad \le
2\alpha_m\bar{\alpha}_m \mathbb{E}_\rho(a(e^{(m)},\tilde{\psi}_{\omega})(\Lambda \|v_\omega\|_\omega - a(u,\tilde{\psi}_{\omega})))\\
&&\qquad \le 2\alpha_m\bar{\alpha}_m \mathbb{E}_\rho(a(e^{(m)},\tilde{\psi}_{\omega})^2)^{1/2}\mathbb{E}_\rho((\Lambda\|v_\omega\|_\omega- a(u,\tilde{\psi}_{\omega}))^2)^{1/2}\\
&&\qquad \le \alpha_m^2 \mathbb{E}_\rho(a(e^{(m)},\tilde{\psi}_{\omega})^2) + \bar{\alpha}_m^2(\Lambda^2\mathbb{E}_\rho(\|v_\omega\|_\omega^2)-2\Lambda
\mathbb{E}_\rho(\|v_\omega\|_\omega\, a(u,\tilde{\psi}_{\omega})) + \mathbb{E}_\rho(a(u,\tilde{\psi}_{\omega})^2))\\
&&\qquad \le \alpha_m^2 \mathbb{E}_\rho(a(e^{(m)},\tilde{\psi}_{\omega})^2) + \bar{\alpha}_m^2(\Lambda^2\|u\|_{A_2}^2+2\Lambda \|u\|_{A_2} \|u\| + \mathbb{E}_\rho(a(u,\tilde{\psi}_{\omega})^2)),
\eea
where we have used that by (\ref{A1})
$$
|\mathbb{E}_\rho(\|v_\omega\|_\omega a(u,\tilde{\psi}_{\omega}))|\le \mathbb{E}_\rho(\|v_\omega\|_\omega)\|u\|\le 
\mathbb{E}_\rho(\|v_\omega\|_\omega^2)^{1/2}\|u\| = \|u\|_{A_2}\|u\|.
$$
After substitution into (\ref{I1}) some terms cancel, and  we arrive at the estimate
\be\label{I2}
\mathbb{E}_\rho(\delta_{m+1}^2|u^{(m)})\le \alpha_m^2 \delta_m^2 + \bar{\alpha}_m^2(\Lambda \|u\|_{A_2}+\|u\|)^2
\ee
for the expectation of the squared error $\delta_{m+1}^2$ conditioned on $u^{(m)}$. Because of the independence assumption, this gives
the recursion for the expected error
\be\label{ERec}
\mathbb{E}(\delta_{m+1}^2) \le \alpha_m^2 \mathbb{E}(\delta_m^2) + \bar{\alpha}_m^2(\Lambda \|u\|_{A_2}+\|u\|)^2,\qquad
m=0,1,\ldots,
\ee
with $\mathbb{E}_\rho(\delta_{0}^2)=\|u\|^2$ since we set $u^{(0)}=0$. Due to the specific choice of $\alpha_m$, 
for the sequence $b_m:=(m+1)\mathbb{E}(\delta_m^2)$ this yields the recursion
$$
b_{m+1} \le \alpha_m b_m + \bar{\alpha}_m(\Lambda\|u\|_{A_2}+\|u\|)^2, \qquad m=0,1,\ldots,\qquad b_0=\|u\|^2,
$$
which implies $b_m\le (\Lambda\|u\|_{A_2}+\|u\|)^2$ uniformly in $m$ (note $\alpha_m+\bar{\alpha}_m=1$). This is equivalent to (\ref{EC2}), and concludes the proof  of Theorem \ref{theo1}.\hfill $\Box$

\smallskip
As in \cite{GrOs2016}, the proof of Theorem \ref{theo1} can be modified to yield an estimate valid for arbitrary $u\in V$.
This results in the following:

\begin{theo}\label{theo2}
If the functions defined in (\ref{A1}) are Bochner-measurable then for arbitrary $u\in V$ the algorithm (\ref{Rec})
satisfies
\be\label{ECV}
\mathbb{E}(\delta_{m}^2)^{1/2} \le 2\left(\|u-h\| +\frac{((\Lambda\|h\|_{{A}_2}+\|h\|)^2+\|u\|^2)^{1/2}}{(m+1)^{1/2}}\right),\qquad m=0,1,\ldots,
\ee
where $h\in {A}_2$ is arbitrary. As a consequence, we have $\mathbb{E}(\delta_{m}^2)\to 0$ if $u$ belongs to the closure of $A_2$ in $V$.
\end{theo}

{\bf Proof}. To see (\ref{ECV}), write
$$
a(e^{(m)},u)-a(e^{(m)},\tilde{\psi}_{\omega_m})\,a(u,\tilde{\psi}_{\omega_m})\le a(e^{(m)},h)-a(e^{(m)},\tilde{\psi}_{\omega_m})\,a(h,\tilde{\psi}_{\omega_m})
+\|u-h\|\|e^{(m)}\|,
$$
and proceed as above for the first term in the right-hand side, using the assumption $h\in {A}_2$. Instead of (\ref{ERec}), this yields
\be\label{ERecV}
\mathbb{E}(\delta_{m+1}^2) \le \alpha_m^2 \mathbb{E}(\delta_m^2) + 2\alpha_m\bar{\alpha}_m \mathbb{E}(\delta_m)\|u-h\|+\bar{\alpha}_m^2((\Lambda\|h\|_{A_2}+\|h\|)^2+\|u\|^2),
\ee
$m=0,1,\ldots$. The rest of the argument leading to (\ref{ECV}) is the same as in the proof of \cite[Theorem 2]{GrOs2016}. \hfill $\Box$

\smallskip
The bounds in Theorem \ref{theo1} and Theorem \ref{theo2} carry over to the stochastic version of orthogonal matching pursuit (OMP), where the recursion (\ref{Rec}) is replaced by
\be \label{OMP}
u^{(m+1)} = P_{W_{m}}u,\qquad r^{(m)}_{\omega_m}=T_{\omega_m}(u-u^{(m)}),\quad m=0,1,\ldots,\quad u^{(0)}=0,
\ee
with $P_{W_{m}}$ denoting the orthogonal projection onto the subspace 
$$
W_{m}:=\mathrm{span}(\{R_{\omega_0}r^{(0)}_{\omega_0},\ldots,R_{\omega_m}r^{(m)}_{\omega_m}\})
$$
in $V$. This is because, for given $u^{(k)}$ and $\omega_k$, $k=0,\ldots,m$, the error of the stochastic OMP algorithm after  the update step satisfies
$$
\|u-u^{(m+1)}\|=\|u-P_{W_{m}}u\|\le \|u-\alpha_m u^{(m)}-\xi_m R_{\omega_m}r^{(m)}_{\omega_m}\|
$$
for any choice of $\xi_m$. Consequently, the estimates for one step of (\ref{Rec}) can be applied, and we obtain the same
recursions for the expectations of the squared error $\mathbb{E}(\|u-P_{W_{m}}u\|^2)$ of stochastic OMP as 
in (\ref{ERec}) and (\ref{ERecV}) for our algorithm (\ref{Rec}). Thus, the bounds in Theorem \ref{theo1} and 
Theorem \ref{theo2} hold for stochastic OMP as well. In practice, stochastic OMP (\ref{OMP}) is expected to converge slightly faster than our algorithm (\ref{Rec}), at the expense of a more costly evaluation of the projections 
$P_{W_{m}}u$ in each step.

Finally, note that the class ${A}_2\subset V$ delicately depends on the choices for $\{V_\omega\}_{\omega\in \Omega}$ and the probability measure $\rho$. It can be made more explicit in some cases which we outline in the next section.

\section{Examples}\label{sec3}

\subsection{Countable $\Omega$ and discrete measures}\label{sec31} The most common situation in which our results can be made more explicit is the case of a discrete measure $\rho$ on a countable $\Omega$. To allow for a direct comparison with the results in \cite{GrOs2016}, set without loss of generality $\Omega=\mathbb{N}$
and denote $\rho_i:=\rho(\{i\})>0$. Then the measurability assumptions for (\ref{A1}) and (\ref{A2}) are irrelevant. Thus,  $u\in {A}_2$ is, according to (\ref{A2}), equivalent to
the existence of $v_i\in V_i$ such that
$$
u=\sum_i \rho_i R_iv_i,\qquad \sum_i \rho_i \|v_i\|^2_i < \infty, 
$$
where $\|\cdot\|_i$ is the norm in $V_i$. Moreover,
$$
\|u\|_{A_2}^2 = \inf_{u=\sum_{i} R_iv_i}\, \sum_i \rho_i \|v_i\|^2_i .
$$
In \cite{GrOs2016}, for any sequence $\gamma:=\{\gamma_i>0\}$ and $0<q\le \infty$, classes $\mathcal{A}_q^\gamma$ were introduced by the requirement that $u\in \mathcal{A}_q^\gamma$ if there are $w_i\in V_i$ such that
$$
u=\sum_i R_iw_i,\qquad \|\{\gamma_i^{-1}\|w_i\|_i \}\|_{\ell_q}< \infty.
$$
The (quasi-)norm on $\mathcal{A}_q^\gamma$ is given by
$$
\|u\|_{\mathcal{A}_q^\gamma} = \inf_{u=\sum_{i} R_iw_i}\, \|\{\gamma_i^{-1}\|w_i\|_i \}\|_{\ell_q}.
$$
In particular, for the sequence $\gamma=\mathbf{1}$ given by $\gamma_i=1$, we simply 
use the notation $\mathcal{A}_q=\mathcal{A}_q^{\mathbf{1}}$.

\begin{lem}\label{lem1}
For any given $\{V,a\}$ and $\{V_i,a_i\}_{i\in\mathbb{N}}$ and any discrete probability measure $\rho$ on $\mathbb{N}$, the following continuous embbedings hold with norm $\le 1$:
$$
\mathcal{A}^\rho_\infty \subset A_2=\mathcal{A}_2^{\sqrt{\rho}} \subset \mathcal{A}_1.
$$
\end{lem}

{\bf Proof}. Since
$$
\|u\|_{\mathcal{A}^{\sqrt{\rho}}_2}^2=\inf_{u=\sum_{i} R_iw_i} \sum_i \rho_i^{-1}\|w_i\|_i^2 =
\inf_{u=\sum_{i} \rho_i R_iv_i} \sum_i \rho_i\|v_i\|_i^2 =\|u\|^2_{A_2},
$$
the equality $A_2=\mathcal{A}_2^{\sqrt{\rho}}$ is obvious. Take any $u$ of the form $u=\sum_i R_i w_i$.  The inequalities
$$
\sum_i \rho_i^{-1}\|w_i\|_i^2 \le \sum_i \rho_i \sup_i (\rho_i^{-1}\|w_i\|_i)^2=(\sup_i \rho_i^{-1}\|w_i\|_i)^2
$$
and 
$$
\sum_i \|w_i\|_i = \sum_i \rho_i^{1/2}(\rho^{-1/2}\|w_i\|_i) \le (\sum_i \rho_i)^{1/2}(\sum_i \rho_i^{-1}\|w_i\|_i^2)^{1/2}
=(\sum_i \rho_i^{-1}\|w_i\|_i^2)^{1/2} 
$$
imply the embeddings $\mathcal{A}^\rho_\infty \subset \mathcal{A}_2^{\sqrt{\rho}}$ and
$\mathcal{A}_2^{\sqrt{\rho}} \subset \mathcal{A}_1$, respectively.\hfill $\Box$

\smallskip
In \cite{GrOs2016}, the condition $u\in \mathcal{A}^\rho_\infty$ was shown to be sufficient for estimates essentially identical with (\ref{ERec}) and
(\ref{ERecV}) to hold, therefore the present paper improves the results from \cite{GrOs2016} (and extends them to uncountable $\Omega$).
On the other hand, as shown in \cite{BCDD2008,GrOs2016} the condition $u\in \mathcal{A}_1$ is sufficient for proving convergence rates similar to (\ref{ERec}) and
(\ref{ERecV}) for the weak greedy version of our algorithm (\ref{Rec}), where the random choice of $\omega_m$ is replaced by a residual-based search for
a $\omega_m\in \Omega$ such that
\be\label{Gr}
\|r_{\omega_m}^{(m)}\|_{\omega_m} \ge \beta \sup_{\omega\in\Omega}\|r_{\omega}^{(m)}\|_{\omega}.
\ee
Here, $\beta\in (0,1]$ is a fixed parameter. In other words, for the specific algorithm (\ref{Rec}) the greedy rule (\ref{Gr}) 
of picking the $\omega_m$ yields the same convergence bound on a larger class of $u$ than any of the stochastic search algorithms. The drawback of greedy algorithms is the cost of
implementing (\ref{Gr}) which typically requires the computation of residuals $r_\omega^{(m)}$ for many $\omega\in \Omega$.

\subsection{Stochastic approximation in RHKS}\label{sec32} Another case where the above theory can be substantiated is the approximation of functions in a reproducing kernel Hilbert space from randomly selected point evaluations. The standard setting \cite{RKHS,Bog} is as follows: Let $\Omega$ be a compact metric space, and let
$K:\,\Omega\times \Omega \to \mathbb{R}$ be a continuous positive-definite kernel. This kernel defines a Hilbert space
$H_K$ with scalar product $(\cdot,\cdot)_K$ whose elements are continuous functions $f:\;\Omega\to \mathbb{R}$ such that
\be\label{RK}
(K_\omega,f)_K = f(\omega) \qquad \forall\;f\in H_K\quad \forall \;\omega\in\Omega.
\ee
Here, $K_\omega\in H_K$ is given by $K_\omega(\eta)=K(\omega,\eta)$, $\eta\in \Omega$. 
Now, choose $V=H_K$ with the scalar product $a(\cdot,\cdot)=(\cdot,\cdot)_K$ and consider the family of one-dimensional subspaces 
$V_\omega\subset V$ spanned by $K_\omega$, $\omega\in \Omega$. In particular, $a_\omega(\cdot,\cdot)$ is the restriction of
$(\cdot,\cdot)_K$ to $V_\omega$, and $R_\omega$ is the natural injection ($\Lambda=1$). With this, we compute
$$
R_\omega T_\omega f = T_\omega f = \frac{(K_\omega,f)_K}{(K_\omega,K_\omega)_K}K_\omega =\frac{f(\omega)}{K(\omega,\omega)}K_{\omega},
$$
where in the last step we have used the reproducing kernel property (\ref{RK}). Thus, our algorithm (\ref{Rec}) turns into an incremental 
approximation process, requiring in each step the evaluation of 
$$
e^{(m)}(\omega_m)=f(\omega_m)-u^{(m)}(\omega_m),
$$
where $\omega_m$ is chosen randomly and
independently from $\Omega$ according to a certain probability distribution $\rho$. This scenario is typical in learning theory \cite{SmZh}, where 
the samples $(\omega_m,y_m)\in \Omega\times \mathbb{R}$, which are drawn according to an (unknown) joint probability distribution $\tilde{\rho}$ on 
$\Omega\times \mathbb{R}$, become incrementally available, and one tries to recover the regression function
$$
f(\omega)=\mathbb{E}_{\tilde{\rho}}(y|\omega).
$$ 
In the "no-noise" case ($\mathbb{E}_{\tilde{\rho}}((y-f(\omega))^2|\omega) = 0$ a.e. on $\Omega$), we would have $y_m=f(\omega_m)$ almost surely,
while the $\omega_m$ are independent samples drawn from $\Omega$  according to the marginal distribution $\rho=\tilde{\rho}_\omega$.

To apply our theory, i.e., to obtain rates for the expectation of the squared error from (\ref{EC2}) and (\ref{ECV}), we need to check (\ref{A1}) and have to examine the condition $u\in A_2$ and the approximability of $u\in H_K$ by elements from $A_2$, respectively. The measurability assumptions for (\ref{A1}) and
(\ref{A2}) follow from the uniform continuity of the kernel which implies the uniform continuity of the function $\omega\to K_{\omega}$, and
the measurability of the function $\omega \to R_\omega v_\omega = c_\omega K_\omega$ for any measurable scalar-valued function $\omega\to c_\omega$.
Thus, $u\in A_2$ if 
$$
u(\eta)=(u,K_\eta)_K=\mathbb{E}_\rho( (c_\omega K_\omega,K_\eta)_K)=\int_\Omega c_{\omega}K(\omega,\eta)\,d\rho_\omega,\qquad \int_\Omega c_\omega^2 \,d\rho_\omega < \infty,
$$
i.e., if $u$ is in the image of $L_2(d\rho)$ under the action of the integral operator $L_K$ with kernel $K$
given by the formula
$$
(L_Kf)(\eta) := \int_\Omega K(\omega,\eta)f(\omega)\,d\rho_\omega.
$$ 
It is well known 
that the operator $L_K$ is also well-defined on $V=H_K$, that it is trace-class positive semi-definite on $H_K$, and that 
$A_2=L_K(L_2(d\rho))=L_K^{1/2}(H_K)$. Thus, our result recovers rates for the noiseless case analogous to those known in online learning 
with kernels for similar approximation algorithms \cite{TaYa2013,LiZh2015,DeBa2016}, where the spaces defined in terms of the spectral decomposition of $L_K$ often serve as smoothness classes.

\subsection{General unit norm dictionaries}\label{sec33} As a third, slightly different but also slightly more general example, let us consider the case when, for a given separable Hilbert space
$V=H$ with scalar product $a(\cdot,\cdot)=(\cdot,\cdot)$, we choose a Borel measure $\rho$ concentrated on the unit sphere $\Omega=S_H=
\{\omega\in H:\;\|\omega\|=1\}$ of $H$. Then, we consider the algorithm (\ref{Rec}) with the family $V_\omega :=\mathrm{span}(\{\omega\})$ of one-dimensional
subspaces of $H$ (again, $a_\omega(\cdot,\cdot)=(\cdot,\cdot)$ on $V_\omega$, $R_\omega$ are the natural injections,
and $\Lambda=1$ ). Since any function of the form $\omega \in S_H \to v_\omega=c_\omega \omega$ is Bochner-measurable if the scalar-valued function  $\omega \to c_\omega$ is measurable, we have
$u\in A_2$ iff
\be\label{A2a}
u=\int_{S_H} c_\omega \omega d\rho, \qquad \int_{S_H} c_\omega^2\, d\rho_\omega < \infty.
\ee
In this case, the proof of (\ref{EC2}) can be carried out directly, using the covariance operator $L:\,H\to H$ given by
\be\label{CoOp}
L v = \mathbb{E}_\rho((v,\omega)\omega) =\int_{S_H} (v,\omega)\omega \,d\rho_\omega, \qquad v\in H.
\ee
This operator is positive semi-definite and trace-class, i.e., there is a complete orthonormal system of eigenfunctions $\psi_k$ of $L$ for the subspace 
$$
\tilde{H}:= H\ominus \mathrm{Ker}(L)
$$
with associated eigenvalues $\mu_k>0$ satisfying $\sum_k \mu_k=1$. The powers $L^s$, $s>0$, are well defined on $H$ and act
as isometries between $\tilde{H}$ and the Hilbert spaces 
$$
H^s_L=L^s(H):=\{v=\sum_k \mu_k^s c_s \psi_k: \; \|v\|_{H^s_L}:= (\sum_k c_k^2)^{1/2}<\infty\}.
$$
The latter serve as smoothness spaces and, as we will see, $u\in H^{1/2}_L$ implies an analog of (\ref{EC2}).
Indeed, since $\omega\in S_H$ we have $\tilde{\psi}_{\omega}=\omega$ in (\ref{A1}) for any $e$.  Taking into account (\ref{CoOp}) the counterpart of (\ref{I1}) reads as follows:
\bea
&&\mathbb{E}_\rho(\delta_{m+1}^2)= \alpha_m^2 (\delta_m^2-\mathbb{E}_\rho((e^{(m)},\omega)^2))\\
&&\qquad\qquad\qquad\qquad +2\alpha_m\bar{\alpha}_m((e^{(m)},u)-\mathbb{E}_\rho((e^{(m)},\omega)(u,\omega)))
+\bar{\alpha}_m^2(\|u\|^2-\mathbb{E}_\rho((u,\omega)^2))\\
&&\qquad\qquad \,= \alpha_m^2 (\delta_m^2-(Le^{(m)},e^{(m)}))+2\alpha_m\bar{\alpha}_m((e^{(m)},u)-(Le^{(m)},u))+\bar{\alpha}_m^2(\|u\|^2-(Lu,u)).
\eea
Assuming $u\in H^{1/2}_L$, i.e., $u=L^{1/2}v$ for some $v\in \tilde{H}\subset H$ with $\|u\|_{H^{1/2}_L}=\|v\|$, we estimate the second term
in the right-hand side by
\bea
2\alpha_m\bar{\alpha}_m((e^{(m)},u)-(Le^{(m)},u))&=&2\alpha_m\bar{\alpha}_m(L^{1/2}e^{(m)},(L^{-1/2}-L^{1/2})u)\\
&\le& 2\alpha_m\bar{\alpha}_m \|L^{1/2}e^{(m)}\|\|(L^{-1/2}-L^{1/2})u\|\\
&\le& \alpha_m^2  (Le^{(m)},e^{(m)}) +\bar{\alpha}_m^2 (\|v\|^2 -2\|u\|^2 + (Lu,u)).
\eea
Substitution and cancellation of several terms yields the following analog of (\ref{I2}):
$$
\mathbb{E}_\rho(\delta_{m+1}^2)\le \alpha_m^2 \delta_m^2 + \bar{\alpha}_m^2 \|u\|_{H^{1/2}_L}^2.
$$
The rest is as in the above proof of Theorem \ref{theo1}. This results in the following estimate with slightly improved constant.

\begin{theo}\label{theo3} In the setting described in this subsection,  the algorithm (\ref{Rec}) converges in expectation 
for arbitrary $u\in H_L^{1/2}$:
\be\label{EC2a}
\mathbb{E}(\delta_{m}^2)\le \frac{\|u\|_{H^{1/2}_L}^2}{m+1},\qquad m=0,1,\ldots .
\ee
The analog of (\ref{ECV}) is 
\be\label{ECVa}
\mathbb{E}(\delta_{m}^2)^{1/2}\le 2(\|u-h\|+\frac{(\|h\|_{H^{1/2}_L}^2+\|u\|^2)^{1/2}}{(m+1)^{1/2}}),\qquad m=0,1,\ldots ,
\ee
valid for any $u\in H$ and $h\in H_L^{1/2}$. Convergence in expectation $\mathbb{E}(\delta_{m}^2)\to 0$ holds for any
$u\in \tilde{H}$.\\
Moreover, the classes $A_2$ and $H_L^{1/2}$ coincide, with equality of norms $\|u\|_{A_2}=\|u\|_{H_L^{1/2}}$ for any $u\in H_L^{1/2}$.
\end{theo}

{\bf Proof}. The estimate (\ref{EC2a}) was already established, the modification leading to
(\ref{ECVa}) is similar to the one in the proof of Theorem \ref{theo2}: Since
\bea
(e^{(m)},u)-(Le^{(m)},u)&=&(e^{(m)},h)-(Le^{(m)},h) + (e^{(m)},(I-L)(u-h))\\
&\le& (e^{(m)},h)-(Le^{(m)},h) + \|e^{(m)}\|\|u-h\|,
\eea
we can proceed for the first term as above, with $u$ replaced by $h\in A_2$,
to arrive at
$$
\mathbb{E}(\delta_{m+1}^2)\le \alpha_m^2 \mathbb{E}(\delta_m^2)+2\alpha_m\bar{\alpha}_m\mathbb{E}(\delta_m)\|u-h\| +\bar{\alpha}_m^2(\|h\|^2_{A_2}+\|u\|^2).
$$
The last term results from a rough estimate of the collection of all terms with forefactor $\bar{\alpha}_m^2$
remaining after substitution, namely
\bea
&&\|h\|_{A_2}^2-2\|h\|^2 +(Lh,h) +\|u\|^2-(Lu,u)=\|h\|_{A_2}^2+\|u\|^2-\|h\|^2 -((I-L)h,h) -(Lu,u)\\
&&\qquad\qquad \le \|h\|_{A_2}^2+\|u\|^2.
\eea
For the rest of the argument,
we again refer to the proof of Theorem 2 b) in \cite{GrOs2016}.

It remains to check that $A_2=H_L^{1/2}$. For $u\in A_2$ satisfying (\ref{A2a}) we can write
$$
\|u\|_{H_L^{1/2}}^2 =\sum_k \frac{(u,\psi_k)^2}{\mu_k} = \sum_k \left(\int_{S_H} c_\omega (\omega,\mu_k^{-1/2}\psi_k)\, d\rho_\omega\right)^2
= \sum_k (c_\omega,f_{k,\omega})^2_{L_2(d\rho)} \le \|c_\omega\|_{L_2(d\rho)}^2 <\infty.
$$
The last step follows because the functions $f_{k,\omega}:=(\omega,\mu_k^{-1/2}\psi_k)$ form an orthonormal system in
$L_2(d\rho)$:
$$
(f_{k,\omega},f_{l,\omega})^2_{L_2(d\rho)} = \int_{S_H} \frac{(\omega,\psi_k))(\omega,\psi_l)}{\mu_k^{1/2}\mu_l^{1/2}}\,d\rho_\omega
=\frac{(L\psi_k,\psi_l)}{\mu_k^{1/2}\mu_l^{1/2}} = \delta_{kl}.
$$
Moreover, for similar reasons any $u\in A_2$ must be orthogonal to $\mathrm{Ker}(L)$, i.e., belongs to $\tilde{H}$ and is thus in the 
closure in $H$ of the orthonormal system $\{\psi_k\}$ of eigenfunctions of $L$. Indeed, if $v\in \mathrm{Ker}(L)$ then we have $(\omega,v)=0$ almost everywhere on $\Omega$ since
$$
\int_{S_H} (\omega,v)^2 \,d\rho_\omega = (Lv,v)=0.
$$
This implies the desired orthogonality
$$
(u,v)=\int_{S_H} c_\omega (\omega,v) \,d\rho_\omega = 0,
$$
and shows $u\in H^{1/2}_L$ and $\|u\|_{H^{1/2}_L}\le \|u\|_{A_2}$ for all $u\in A_2$. 
 
Now, take $u\in H^{1/2}_L$, i.e.,
$$
u=\sum_k c_k \psi_k, \qquad \|u\|_{H^{1/2}_L}^2 = \sum_k \mu_k^{-1}c_k^2 < \infty.
$$
We will check that (\ref{A2a}) holds with $c_\omega = \sum_k \mu_k^{-1/2}c_k f_{k,\omega}$, which immediately implies
$u\in A_2$ and the opposite inequality $\|u\|_{A_2}\le \|u\|_{H^{1/2}_L}$. This is done by verifying that the moments
$(u,\psi_l)$ coincide for both representations of $u$: On the one hand, we have $(u,\psi_l)=c_l$, on the other hand,
we have
$$
\left(\int_{S_H} c_\omega\omega\,d\rho_\omega,\psi_l\right)=\int_{S_H} \left(\sum_k \mu_k^{-1/2}c_k f_{k,\omega}\right) (\omega,\psi_l) \,d\rho_\omega
=\sum_k (\mu_l/\mu_k)^{1/2} c_k (f_{k,\omega},f_{l,\omega})_{L_2(d\rho)} = c_l
$$
by the orthonormality of the system $\{f_{k,\omega}\}$ in $L_2(d\rho)$.\hfill $\Box$

\subsection{Collective approximation}\label{sec34} To demonstrate the versatility of the abstract scheme developed in Section \ref{sec2}, we consider a problem raised in \cite{AJOP2017}: Given 
an $n$-dimensional subspace $V_n$ of a Hilbert space $H$ and a dictionary $D$ of unit norm elements in $H$
(the condition $D\subset S_H$ is silently kept througout this subsection),  
construct, by incrementally selecting dictionary elements $\omega_0,\omega_1,\ldots$, subspaces $W_{m-1}=\mathrm{span}\{\omega_0,\ldots,\omega_{m-1}\}$ which approximate $V_n$ well, i.e., for which estimates for the approximation quantities
$$
\sigma_m=\sup_{v\in V_n:\,\|v\|=1} \inf_{w\in W_{m-1}} \|v-w\|_H = \sup_{v\in V_n:\,\|v\|=1} \|v-P_{W_{m-1}}v\|_H
$$
hold.
The collective OMP algorithm proposed in \cite{AJOP2017} uses greedy selection of $\omega_m\in D$ based on computations involving the ortho-projections $P_{W_{m-1}}$ onto $W_{m-1}$ which become more costly for larger $m$. It comes with a convergence rate for the quantity
$$
\epsilon_m({\Phi})= \left(\sum_{i=1}^n \|\phi_i -P_{W_{m-1}}\phi_i\|_H^2\right)^{1/2} = \|\Phi-P_{W_{m-1}}\Phi\|_{H^n} \ge \sigma_m,\qquad m=1,2,\ldots,
$$
where ${\Phi}=(\phi_1,\ldots,\phi_n)$ is a given orthonormal  basis in $V_n$.

We apply our results and design algorithms avoiding the projections $P_{W_{m-1}}$ while still guaranteeing similar convergence rates. To set the scene, identify $V$ with $H^n$ equipped with the usual scalar product
$$
a(\mathbf{u},\mathbf{v}):=\sum_{i=1}^n (u_i,v_i),\qquad \mathbf{u},\mathbf{v}\in V \quad (\;\mathbf{u}=(u_1,\ldots,u_n)\;).
$$
Let $\Omega=D$, and consider the family
$$
V_\omega := \{\mathbf{v}_\omega=\mathbf{c}\omega:\; \mathbf{c}\in\mathbb{R}^n\},\qquad \omega\in \Omega,
$$
of $n$-dimensional subspaces of $V$ (again, $R_\omega$ are the natural injections, $\Lambda=1$). The problem we want to solve is $\mathbf{u}=\Phi$ or, in variational form,
$$
a(\mathbf{u},\mathbf{v})=a(\Phi,\mathbf{v})\qquad \forall\; \mathbf{v}\in \mathbf{u}.
$$
With this, we have
$$
R_\omega T_\omega \mathbf{v} = T_\omega \mathbf{v} =a (\mathbf{v},\omega)\omega,
$$
where $a(\mathbf{v},\omega):=((v_1,\omega),\ldots,(v_n,\omega))\in \mathbb{R}^n$. 

Independently of the method of choosing $\omega_m$ (randomly or greedy), our algorithm (\ref{Rec})
$$
\mathbf{u}^{(m+1)} = \alpha_m \mathbf{u}^{(m)} + \xi_m r^{(m)}_{\omega_{m}},\qquad
r^{(m)}_{\omega_{m}}=T_{\omega_{m}}\mathbf{e}^{(m)}= (\Phi-\mathbf{u}^{(m)},\omega_m)\omega_{m},
\qquad m=0,1,\ldots,
$$
when started with $\mathbf{u}^{(0)}=\mathbf{0}$,
produces a sequence of $\mathbf{u}^{(m)}$ whose components belong to $W_{m-1}$ if $m>0$. Thus, we have upper estimates
$$
\epsilon_m({\Phi})\le \delta_m:=\|\Phi-\mathbf{u}^{(m)}\|,\qquad m=1,2,\ldots.
$$

If we choose the $\omega_m$, $m=0,1,\ldots$, randomly and independently according to a Borel measure $\rho$ on $S_H$  with support on $D$, then Theorems \ref{theo1} and \ref{theo2} are applicable, and they imply rates (in expectation)
for $\Phi\in A_2$ and general $\Phi$ in terms of its approximability by elements $\mathbf{h}\in A_2$. 
Moreover, it is easy to see that the proof of Theorem \ref{theo3} remains valid  if the application of the operators 
$L$ and $L^s$, respectively, which are defined on $H$ and depend on $\rho$, is extended componentwise to
$V=H^n$. This way, we obtain the estimate
\be\label{EC2b}
\sigma_m^2\le \epsilon_m^2\le \mathbb{E}(\delta_{m}^2)\le \frac{\|\Phi\|_{A_2}^2}{m+1},\qquad m=1,2,\ldots ,
\ee
if $\Phi\in A_2=(H^{1/2}_L)^n$ with norm in $A_2$ defined as
$$
\|\mathbf{v}\|_{A_2}^2 =\sum_{i=1}^n \|u_i\|^2_{H^{1/2}_L}.
$$
The counterpart of (\ref{ECVa}) holds, too: If $\Phi\in H^n$ then for arbitrary $\Psi\in A_2$ we have
\be\label{ECVb}
 \mathbb{E}(\delta_{m}^2)^{1/2}\le 2(\|\Phi-\Psi\| + \frac{(\|\Psi\|_{A_2}^2+\|\Phi\|^2)^{1/2}}{(m+1)^{1/2}}),
\qquad m=1,2,\ldots.
\ee
 
 These estimates for the expected error decay of our randomized algorithm are qualitatively the same as for the more expensive collective OMP algorithm with weak greedy selection of the $\omega_m$
 proposed in \cite{AJOP2017}. However, the class $A_2$ is smaller  then the class $\mathcal{A}_1(D)$ appearing in the convergence theory in  \cite{AJOP2017}, and depends on the choice for $\rho$. 
 
 The weak greedy version of our algorithm was already analyzed in  \cite{GrOs2016} by generalizing earlier results from \cite{BCDD2008}. For completeness, we repeat it here in the setting and notation of Section \ref{sec2}.
 
Define the class $A_1$ as the set of all $u\in V$ for which a representation of the form
 \be\label{A11}
 u=\sum_j R_{\omega^j}v_{\omega^j}, \qquad \sum_j\|v_{\omega^j}\|_{\omega^j} <\infty,\quad \omega^j\in \Omega,
 \ee
 holds, and set
 $$
 \|u\|_{A_1}:= \inf_{u=\sum_j R_{\omega^j}v_{\omega^j}}\; \sum_j\|v_{\omega^j}\|_{\omega^j} .
 $$
For countable $\Omega$, $A_1$ coincides with the class $\mathcal{A}_1$ defined before.

\begin{theo}\label{theo4}
If $u\in A_1$, the algorithm (\ref{Rec}) with $\omega_m$ chosen according to the weak greedy rule (\ref{Gr}) 
possesses the error bound
\be\label{CG1}
\delta_m^2 \le \frac{2((\Lambda/\beta)^2\|u\|_{A_1}^2 +\|u\|^2)}{m+1},\qquad m=0,1,\ldots.
\ee
\end{theo}
{\bf Proof}.
The proof is almost identical to that of Theorem \ref{theo1}. Indeed, using
 (\ref{Gr}) in (\ref{I0}), we have
 $$
 \frac{\Lambda}{\beta} a(e^{(m)},\tilde{\psi}_{\omega_m}) \ge\frac1{\beta} a_{\omega_m}(r_{\omega_m}^{(m)},\frac{r_{\omega_m}^{(m)}}{\|r_{\omega_m}^{(m)}\|_{\omega_m}})\ge a_{\omega}(r_{\omega}^{(m)},\frac{r_{\omega}^{(m)}}{\|r_{\omega}^{(m)}\|_{\omega}})\ge \frac{a(e^{(m)},R_\omega v_\omega)}{\|v_\omega\|_\omega},
$$
for any $\omega\in \Omega$ (as before $\tilde{\psi}_{\omega_m}$ is defined in (\ref{A1}) with $e=e^{(m)}$). Thus, representing $u\in A_1$ as in (\ref{A11}), we arrive at
$$
a(e^{(m)},u)= \sum_j (a(e^{(m)},R_{\omega^j} v_{\omega^j}) \le \frac{\Lambda a(e^{(m)},\tilde{\psi}_{\omega_m}) }{\beta} \sum_j \|v_{\omega^j}\|_{\omega^j},
$$
and, after taking the infimum over all such representations of $u$, we get
$$
a(e^{(m)},u) \le \frac{\Lambda \|u\|_{A_1} }{\beta}a(e^{(m)},\tilde{\psi}_{\omega_m}).
$$
For the corresponding term of the error representation for $\delta_{m+1}^2$, this yields
\bea
&& 2\alpha_m\bar{\alpha}_m(a(e^{(m)},u)-a(e^{(m)},\tilde{\psi}_{\omega_m})a(u,\tilde{\psi}_{\omega_m})\\
&& \qquad \le 2\alpha_m\bar{\alpha}_m a(e^{(m)},\tilde{\psi}_{\omega_m})((\Lambda/\beta) \|u\|_{A_1}-a(u,\tilde{\psi}_{\omega_m}))\\
&& \qquad \le \alpha_m^2 a(e^{(m)},\tilde{\psi}_{\omega_m})^2+\bar{\alpha}_m^2((\Lambda/\beta)\|u\|_{A_1}-a(u,\tilde{\psi}_{\omega_m}))^2,
\eea
and after substitution and cancellation of terms we have
\bea
\delta_{m+1}^2&\le& \alpha_m^2 \delta_m^2 +\bar{\alpha}_m^2
(((\Lambda/\beta) \|u\|_{A_1}-a(u,\tilde{\psi}_{\omega_m}))^2+\|u\|^2-a(u,\tilde{\psi}_{\omega_m})^2)\\
&\le& \alpha_m^2 \delta_m^2 +2\bar{\alpha}_m^2
((\Lambda/\beta)^2\|u\|_{A_1}^2+\|u\|^2).
\eea
The rest is as before. \hfill $\Box$

\section{Concluding remarks} \label{sec4}
We conclude with three further remarks.

{\bf Remark 1.} In the generality considered here, the obtained convergence rates for the expectation of the squared error $\delta_m^2$ of the algorithm (\ref{Rec}) for $u\in A_2$ cannot be improved without additional assumptions on $\rho$ or $u$. 
To see this, consider the case of a discrete measure $\rho$ concentrated on a complete orthonormal system $\{e_j\}\subset S_H$ in a Hilbert space $V=H$ with scalar product $a(\cdot,\cdot)=(\cdot,\cdot)$, and denote $\rho_j=\rho(\{e_j\}) >0$, $j\in \Omega=\mathbb{N}$.  This is within the setting of Section \ref{sec33}. Obviously, we have 
$$
Lv=\sum_j \rho_j(v,e_j) e_j,\qquad v\in H\quad (\psi_j=e_j, \;\mu_j=\rho_j,\; j\in \mathbb{N}),
$$
and $\mathrm{Ker}(L)=\{0\}$. In other words, $\tilde{H}=H$, and
$$
H^s_L = \{u=\sum_j \rho_j^{s} c_j e_j:\quad \|u\|_{H^s_L}^2:= \sum_j c_j^2<\infty\},\qquad  s\in \mathbb{R}.
$$
Following the reasoning in Remark 5 in \cite{GrOs2016}, for any algorithm that produces 
the iterates $u^{(m)}$ as linear combinations of at most $m$ elements $e_j$ drawn randomly and independently according to $\rho$, we then have  the lower estimate
$$
\mathbf{E}(\|u-u^{(m)}\|^2)\ge \sum_j (u,e_j)^2(1-\rho_j)^m.
$$
We mention as a side note that this lower bound is achieved for the stochastic OMP method (\ref{OMP}). The condition $u \in H^{r}_L$ is for $r>0$ equivalent to
$(u,e_j)=\rho_j^r (v,e_j)$, $j\in \Omega$, for some $v\in H$.
Thus, the worst case behavior of the expected squared error of any such algorithm for recovering $u\in H^{r}_L$ is characterized by 
$$
\epsilon_{m,r}:=\sup_{0\neq u\in H^{r}_L} \frac{\mathbf{E}(\|u-u^{(m)}\|^2)}{\|u\|^2_{H^r_L}}\ge \sup_{0\neq v\in H}
\frac{\sum_j (v,e_j)^2\rho_j^{2r}(1-\rho_j)^m}{\sum_j (v,e_j)^2}=\sup_j \rho_j^{2r}(1-\rho_j)^m.
$$
Since the function $f(t)=t^{2r}(1-t)^m$ takes its maximum for $t\in [0,1]$ at $t_0=(m+2r)^{-1}$, we see that in general no rate better than $\mathrm{O}(m^{-2r})$ can be expected on the class $H^{r}_L$. If we take
$r=1/2$, we see that Theorem \ref{theo3} provides an optimal result, in the sense that the upper limit of $m^{2r}\epsilon_{m,r}$ for $m\to \infty$ is finite \emph{and} strictly positive for any $\rho$.

However, with other assumptions on $u$ or on the spectral properties of $L$ (as it is custom in learning with kernel methods \cite{DeBa2016,LiZh2015}), one may expect better results.

\smallskip\noindent
{\bf Remark 2.} The choice for the parameters $\alpha_m$ and $\xi_m$ in the algorithm (\ref{Rec}) is appropriate if  the evaluations in  (\ref{Rec}) and (\ref{Xi}) (in particular, the functional evaluation $F(u^{(m)}$)) are exact. If one attempts to analyze the same algorithm with, e.g., an independent additive noise term $\varepsilon_m$ the update formula  (\ref{Rec}) (in addition to independence, assume $\mathbb{E}(\varepsilon_m)=0$, and $\sigma^2:=\mathbb{E}(\|\varepsilon_m\|^2)=\mathrm{const.}>0$)  then, in the formulas for $\delta_{m+1}^2$ and subsequently in (\ref{ERec}),
an additional term $\sigma^2$ appears in the right-hand side, i.e.,
$$
\mathbf{E}(\delta_{m+1}^2)\le \alpha_m^2 \mathbf{E}(\delta_{m}^2)+ \bar{\alpha}_m^2(\Lambda \|u\|_{A_2}+\|u\|)^2+\sigma^2,\qquad m=0,1,\ldots.
$$
Now, the term $\sigma^2$ renders any attempt of proving  $\mathbf{E}(\delta_m^2)\to 0$ meaningless.  At the $m$-th step of the recursion, an additional term of the order $(m+1) \sigma^2/3$ would appear in the final estimate for $\mathbf{E}(\delta_m^2)$ which is subdominant only for small $\sigma^2$ and at the initial stages of the iteration. A crude calculation shows that under these assumptions the best possible bound for the expectation of the squared error is
$$
\mathbf{E}(\delta_m^2) \approx \frac{\sigma}{(\Lambda\|u\|_{A_2}+\|u\|)} \quad \mbox{if } \;   m\approx \frac{\Lambda \|u\|_{A_2}+\|u\|}{\sigma}.
$$
This says that, on average, the squared error $\delta_m^2$ cannot be approximated better than the standard deviation of the additive noise $\varepsilon_m$ relative to the size of $u$ which is unsatisfactory. A possible repair is to give up the minimization requirement for $\xi_m$, and  to execute (\ref{Rec}) with some suitably chosen sequence $\xi_m\to 0$.  This is well understood for kernel methods in online learning, and represents a significant difference between the noisy and noiseless case.

\smallskip\noindent
{\bf Remark 3.} \smallskip
The right-hand sides in the estimates (\ref{ECV}) and (\ref{ECVa}) in Theorem \ref{theo2} and Theorem \ref{theo3}
have the form of a $K$-functional for the pairs $(V,{A}_2)$ and $(H,H_L^{1/2})$, respectively, This implies that
rates of the form 
$$
\mathbb{E}(\delta_{m}^2)=\mathrm{O}(m^{-\theta}),\qquad m\to \infty,
$$
with exponent $\theta\in (0,1)$ hold for spaces obtained by real interpolation. 
E.g., in the setting of Theorem \ref{theo3}, we obtain
\be\label{ECVr}
 \mathbb{E}(\delta_{m}^2)\le C  (m+1)^{-2s} \|u\| _{H_L^s}^2,
\qquad m=0,1,\ldots,
\ee
valid for all $u\in H^s_L$ and $0<s<1/2$ with  a certain fixed constant $C$. Indeed, $u\in H_L^s$ can be represented
as
$$
u=\sum_k c_k\mu_k^s \psi_k, \qquad \|u\|_{H^s_L}^2 =\sum_k c_k^2,
$$
and setting
$$
h=\sum_{k:\,\mu_k\ge (m+1)^{-1}} c_k\mu_k^s \psi_k,
$$
we have
\bea
\|h\|_{H^{1/2}_L}^2 &=& \sum_{k:\,\mu_k\ge (m+1)^{-1}} c_k^2\mu_k^{2s-1}\le (m+1)^{1-2s}\|u\|_{H^s_L}^2,\\
\|u-h\|^2 &=& \sum_{k:\,\mu_k < (m+1)^{-1}} c_k^2\mu_k^{2s}\le (m+1)^{-2s}\|u\|_{H^s_L}^2,\\
\|u\|^2 &=& \sum_{k} c_k^2\mu_k^{2s}\le \|u\|_{H^s_L}^2.
\eea
Thus, after substitution into (\ref{ECVa}), we get (\ref{ECVr}) with $C=(2(1+\sqrt{2}))^2<24$.

 \section*{Acknowledgement}
M. Griebel was partially supported by the project {\em EXAHD} of the 
DFG priority program 1648 {\em Software for Exascale Computing" (SPPEXA)}
and by the Sonderforschungsbereich 1060 {\em The Mathematics of Emergent Effects}
funded by the Deutsche Forschungsgemeinschaft.
This paper was written while P. Oswald held a Bonn Research Chair sponsored by the
Hausdorff Center for Mathematics at the University of Bonn funded by the Deutsche Forschungsgemeinschaft. He is grateful for this support.

\bibliographystyle{amsplain}

\end{document}